\documentclass[oneside]{amsart}
\usepackage{amsfonts}
\usepackage{amssymb}
\usepackage{amsxtra}
\usepackage{amstext}
\usepackage[english]{babel}

\newcommand{\R}{\mathbb{R}}

\newcommand{\Z}{\mathbb{Z}}
\newcommand{\Q}{\mathbb{Q}}
\newcommand{\pp}{\mathbb{P}}

\newcommand{\kF}{\mathcal{F}}

\newtheorem {theo} {Theorem} [section]

\newtheorem{lem}[theo]{Lemma}
\newtheorem{prop}[theo]{Proposition}

\newtheorem{Hypo}[theo]{Hypothesis}

\title[Internal DLA generated by cookie random walks on $\Z$]
      {Internal DLA generated by cookie random walks on $\Z$}

\author{Olivier RAIMOND}
\address{Laboratoire Modal'X, Universit\'e Paris Ouest Nanterre La D\'efense, B\^atiment G, 200 avenue de la R\'epublique 92000 Nanterre, France.}
\email{olivier.raimond@u-paris10.fr}
\author{Bruno SCHAPIRA}
\address{D\'epartement de Math\'ematiques, B\^at. 425, Universit\'e Paris-Sud 11, F-91405 Orsay, cedex, France. }
\email{bruno.schapira@math.u-psud.fr}

\begin{document}

\begin{abstract}
We prove a law of large numbers for the right boundary in the model of internal DLA generated by cookie random walks in dimension one. The proof is based on stochastic algorithms techniques. 
\end{abstract} 

\maketitle

\section{Introduction}

The model of Internal DLA (internal diffusion limited aggregation) on $\Z^d$ is a basic example of cluster growth model, 
which has the Euclidean ball as limiting shape, and for which we have now very good control (presumably sharp) 
on the size of the fluctuations around this limiting shape.

It was introduced, on a mathematical level at least, by Diaconis and Fulton \cite{DF} in dimension one. 
In this case the cluster at any time is an interval, and the analysis of its extremities can be transposed in terms of Polya's urn, 
for which we know that the law of large numbers and the central limit theorem hold. In particular fluctuations around the mean value are diffusive.

In dimension larger than or equal to $2$ it has been first studied by Bramson, Griffeath and Lawler \cite{BGL} who proved that the 
limiting shape is the Euclidean ball and that fluctuations are sublinear. Two years later Lawler \cite{L} showed that fluctuations are subdiffusive 
(with an upper bound in $n^{1/3}$). Finally very recently Jerison, Levine and Sheffield \cite{JLS, JLS2, JLS3}, and independently 
Asselah and Gaudilli\`ere \cite{AG, AG2} have obtained a bound in $\ln n$ in dimension $2$ and in $\sqrt {\ln n}$ in dimension larger than or equal to $3$, 
which agrees with computer simulations.

In its simplest form the model is defined as follows. First at time $0$, the cluster $A(0)$ is empty. Next if the cluster at time $n$ 
is some subset $A(n)$ of $\Z^d$, run a simple random walk and stop it when it first exit $A(n)$. Then $A(n+1)$ is the union of $A(n)$ 
and the site where this random walk stopped. A lot of relevant variations on this model can be imagined. For instance Levine and Peres \cite{LP} 
assume that random walks can start from different sources (not only the origin). 
Another natural modification is to change the underlying graph. For instance Blach\`ere and  Brofferio \cite{BB} 
consider graphs with exponential growth and Shellef \cite{S} looked at the case of a random graph, 
namely a supercritical cluster of percolation.

In another direction it is also interesting to see how the analysis can be transposed when we change the simple random walk 
by some other random process. For instance Enriquez, Lucas and Simenhaus \cite{ELS} study the case of random walks in random environment 
in dimension one and the purpose of this paper is to study the case of cookies random walks (also called multi-excited random walks) in dimension one as well. 
These last processes are simple examples of self-interacting random walks, which have been first introduced and studied by 
Benjamini and Wilson \cite{BW} and Zerner \cite{Zer1}. Actually our purpose is to introduce some new technique for the study of Internal DLA 
in dimension one and then look at a specific example. The case of excited random walks seemed interesting to us, because it contains many different situations 
when we let the parameters of the model vary. 
But most of our analysis carry out for rather general processes, 
and in particular when the random walk is recurrent and satisfies a central limit theorem (see Theorem \ref{generalLLN} and Proposition 
\ref{prophnh} below). For cookies random walks, this has been obtained by Dolgopyat \cite{D} for some values of the parameters.

Our tools are stochastic algorithms techniques which extend those used by Diaconis and Fulton in their original paper.  
Our main result is a law of large numbers for the extremities of the cluster. We first present a general result in Subsection \ref{subgene} and then 
we consider in more details the case of cookies random walks in Subsection \ref{subcookie}. The proofs are postponed to the remaining sections.

\section{Model and results} 
\subsection{A general law of large numbers} 
\label{subgene}
The internal DLA process is a sequence of growing clusters (intervals in dimension one) 
which are constructed via a collection $(X^{(n)}_k,k\ge 0)_{n\ge 1}$, of independent random processes on $\Z$. 
We assume here that these processes start from $0$ and only jump to the nearest neighbor,
which means that for all $n$, $X^{(n)}_0=0$, and for all $k$, $X^{(n)}_{k+1}-X^{(n)}_k = \pm 1$. 
We then define the right boundaries $(d_n,n\ge 0)$ of the intervals as follows. First $d_0=0$, and  
for all $n\ge 0$, $d_{n+1}=d_n+1$, if $X^{(n+1)}$ hits the site $d_n+1$ before the site $d_n-n-1$. 
Otherwise $d_{n+1}=d_n$. (We also implicitly assume that our processes visit a.s. an infinite number of sites, 
which assures that $d_n$ is well defined for all $n\ge 0$.) Now for $n\ge 0$, set
$$x_n:=\frac{d_n+1}{n+2},$$ 
and for $x\in \Z$, set  
$$T^{(n)}_x=\inf\{k\ge 0 \ :\ X^{(n)}_k = x\}.$$
Define also the function $h_n:[0,1]\to [0,1]$, by 
$$h_n(x)=\pp\left[T^{(n)}_{[(n+2)x]}<T^{(n)}_{[(n+2)(x-1)]}\right],$$
where for any $u\in \R$, $[u]$ denotes the integer with the smallest absolute value, at distance strictly less than $1$ from $u$.  
Note that for all $n\ge 0$, $h_n$ is nonincreasing, $h_n(0)=1$ and $h_n(1)=0$.  
The next result is a general law of large numbers for $d_n$. 
\begin{theo} 
\label{generalLLN}
Assume that there exists a function $h$, such that for any $x\in [0,1]$, $h_n(x)$ converges to $h(x)$, as $n\to \infty$. 
Then there exists $p\in [0,1]$, such that   
$x_n\to p$, as $n\to \infty$, almost surely. Moreover $p$ is characterized by the fact that $h(p-)\ge p\ge h(p+)$.  
\end{theo} 
\noindent A typical situation when this theorem applies is when the $X^{(n)}$ are i.i.d. and transient. 
Indeed in this case, it is immediate that $h_n(x)\to p:=\pp[X^{(n)}\to +\infty]$, for all $x\in (0,1)$. 
Another typical situation when the theorem applies well, in the recurrent regime, is when   
there is an invariance principle: for instance if for some $c>0$, the process defined for $t\ge 0$ by  
$$Y^{(n)}_t:= \frac{X^{(n)}_{[n^ct]}}{n},$$
converges in law in the Skorokhod space to some process $(Y_t,t\ge 0)$, as $n\to \infty$. 
In this case, for $x\in \R$, set 
\begin{eqnarray*}
T_x:=\left\{ \begin{array}{ll} 
              \inf\{t\ge 0\ :\ Y_t\ge x\} & \textrm{if }x\ge 0\\
               \inf\{t\ge 0\ :\ Y_t\le x\} & \textrm{if }x<0,\\  
             \end{array}
\right.
\end{eqnarray*}
A basic observation is the following:  
\begin{prop} 
\label{prophnh}
Assume that there exists $c>0$, such that $(Y^{(n)}_t,t\ge 0)$ converges in law to some process $(Y_t,t\ge 0)$. Assume that $T_1$ and $T_{-1}$ are a.s. finite. 
Assume also that $h$ defined for all $x\in [0,1]$, by $h(x): = \pp[T_x<T_{x-1}]$, is continuous. Then, for any $x\in [0,1]$, $h_n(x)\to h(x)$, as $n\to \infty$.   
\end{prop}
\noindent Note that if $Y$ is a continuous semimartingale, with a nonzero martingale part, then the function $h$ defined in the proposition is automatically continuous.

\subsection{Excited random walks}  
\label{subcookie}
We look now more precisely at the case of excited (also called multi-excited or cookie) random walks, which is a typical nontrivial example 
to which our general result apply, at least in most situations. However we stress that for some values of a parameter, Theorem \ref{generalLLN} does not apply 
and this leads to interesting unsolved questions.

We start by recalling the definition of excited random walks and we will then review some known results which are relevant here.

We define an environment as the set of two sequences $(p_{i,+},i\ge 1)$ and $(p_{i,-},i\ge 1)$ of positive reals in $(0,1)$. 
Then the law of an excited random walk $(X_n,n\ge 0)$ in this environment is defined as follows. First $X_0=0$. 
Then, if $(\kF_n,n\ge 0)$ is the natural filtration of $X$, we have 
\begin{eqnarray*} 
\pp(X_{n+1}=X_n + 1 \mid \kF_n) = \left\{ \begin{array}{ll} 
1/2 & \textrm{if }X_n=0 \\
p_{i,+} & \textrm{if }X_n>0 \textrm{ and }L_n(X_n)=i\\
p_{i,-} & \textrm{if }X_n<0 \textrm{ and }L_n(X_n)=i,
\end{array} 
\right. 
\end{eqnarray*} 
where for all $x$ and $n\ge 1$, $L_n(x)=\#\{k\le n \ :\ X_k=x\}$ is the local time of the walk in $x$ at time $n$. 
Moreover $X$ always jumps to a nearest neighbor, so for all $n$, $\pp(X_{n+1}=X_n + 1 \mid \kF_n)= 1 - \pp(X_{n+1}=X_n - 1 \mid \kF_n)$. 
This definition slightly differs from usual ones, where $p_{i,+}=p_{i,-}$ and where $0$ does not play a special role, 
but technically this does not change anything. The reason why we use this definition is that it gives rise to slightly richer situations. 
For instance the process can be transient only to the right.

Now we need some assumption on the environment. We say that a sequence of reals $(x_i,i\ge 1)$ is nonnegative, resp. nonpositive, if $x_i\ge 0$, resp. $x_i\le 0$, for all $i\ge 1$. We will use the following hypothesis:

\begin{Hypo}
\label{cookiespositifs}
Each of the sequences $(p_{i,+}-1/2,i\ge 1)$ and $(p_{i,-} -1/2,i\ge 1)$ is either nonnegative or nonpositive. 
\end{Hypo} 
\noindent Under this hypothesis 
$$\alpha:=\sum_{i=1}^\infty (2p_{i,+}-1),$$ 
and 
$$\beta:=-\sum_{i=1}^\infty(2p_{i,-}-1),$$
are well defined.

\noindent Consider also the process $Y^{(n)}$ defined for each $n\ge 1$ and for all $t\ge 0$ by
\begin{eqnarray}
\label{yn}
Y^{(n)}_t:=\frac{X_{[nt]}}{\sqrt{n}}.
\end{eqnarray}

\noindent It was proved by Zerner \cite{Zer1} that under Hypothesis \ref{cookiespositifs}, $X$ is recurrent if, and only if, $\alpha \le 1$ and $\beta \le 1$. 
Note that a recurrence criterion has also been proved without assuming Hypothesis \ref{cookiespositifs}, but with the other hypothesis that  $p_{i,\pm}=1/2$ for $i$ large enough (see \cite{KZ}). 
Later on Dolgopyat proved the following central limit theorem:\footnote{Actually Zerner and Dolgopyat assumed that $p_{i,+}=p_{i,-}$ for all $i$, but their proofs extend immediately to the present setting.}

\begin{theo}[Dolgopyat \cite{D}]
\label{Dolgo} 
Assume Hypothesis \ref{cookiespositifs} holds. Assume also that $\alpha<1$ and $\beta<1$. Then the process $(Y^{(n)}_t,t\ge 0)$ converges in the 
Skorokhod space toward the $(\alpha,\beta)$-perturbed Brownian motion $(Y_t,t\ge 0)$ defined by the equation: 
\begin{eqnarray}
\label{PBM} 
Y_t=B_t+\alpha\sup_{s\le t} Y_s + \beta\inf_{s\le t} Y_s,
\end{eqnarray}
where $B$ is a Brownian motion. 
\end{theo}

Perturbed Brownian motions are continuous processes which were introduced by Le Gall and Yor \cite{LGY}.  
Uniqueness in law for the equation \eqref{PBM} was later proved 
by Perman and Werner \cite{PW} for any $\alpha<1$ and $\beta<1$, 
and pathwise uniqueness has been proved by Chaumont and Doney \cite{CD}.

We note also that both Zerner's recurrence criterion and Theorem \ref{Dolgo} were proved in the more general case when the $p_{i,\pm}$ may be random. 
Here we will only consider deterministic environments.

If $\alpha<1$ and $\beta<1$, denote by $\Q^{(\alpha,\beta)}$ the law of the $(\alpha,\beta)$-perturbed Brownian motion starting from $0$ 
and by $h_{(\alpha,\beta)}$ the function defined for $x\in[0,1]$ by 
$$h_{(\alpha,\beta)}(x)=\Q^{(\alpha,\beta)}\left[T_x<T_{x-1}\right],$$
where for any $x$, $T_x$ denotes the hitting time of $x$. 
Observe that $h_{(\alpha,\beta)}$ is continuous and decreasing, thus has a unique fixed point $p(\alpha,\beta)$. Since $h_{(\alpha,\beta)}(1)=0=1-h_{(\alpha,\beta)}(0)$, we have $p(\alpha,\beta)\in (0,1)$.  
Actually $h_{(\alpha,\beta)}(x)=P(Z(1-\alpha,1-\beta)\le 1-x)$, where $Z(1-\alpha,1-\beta)$ denotes a Beta random variable with parameters $1-\alpha$ and $1-\beta$ (see Proposition 4 (iii) in \cite{PW}). So another formula for $h_{(\alpha,\beta)}$ is given by 
\begin{eqnarray} 
\label{formuleh}
h_{(\alpha,\beta)}(x) = \int_0^{1-x} \frac{t^{-\alpha}(1-t)^{-\beta}}{B(1-\alpha,1-\beta)}\, dt,
\end{eqnarray}
with $B(\cdot,\cdot)$ the Beta function. Note that $h_{(\alpha,\beta)}(x)=1-h_{(\beta,\alpha)}(1-x)$, which implies that $p(\alpha,\beta)=1-p(\beta,\alpha)$. 
%In particular, since $t^{-\alpha}(1-t)^{\beta}$ is a decreasing function of $t$, we have $p>1/2$, as soon as $\alpha>0$ or $\beta>0$. 
Formula \eqref{formuleh} shows that $p(\alpha,\beta)$ is a continuous function of $(\alpha,\beta) \in (-\infty,1)^2$, and that $p(\alpha,\beta)\to 1$, 
if $\alpha \to 1$ and $\beta<1$ remains fixed (resp. $p(\alpha,\beta)\to 0$, if $\beta\to 1$ and $\alpha<1$ remains fixed). 
Moreover $p$ is monotone in $(\alpha,\beta)$, in the sense that if $\alpha \le \alpha'$ and $\beta \ge \beta'$, then $p(\alpha,\beta)\le p(\alpha',\beta')$.  

\smallskip

\noindent Now we consider the model of internal DLA when $(X^{(n)},n\ge 0)$ is sequence of i.i.d. excited random walks with the same law as $X$ above. 
The next result gives then a law of large numbers for $d_n$ in different situations: 
\begin{theo} 
\label{LLN}
Assume that Hypothesis \ref{cookiespositifs} holds. If $(\alpha,\beta)\neq (1,1)$, 
then there exists $p\in[0,1]$, such that $x_n\to p$ almost surely. 
\begin{itemize}
\item If moreover $\alpha<1$ or $\beta<1$, then $p$ only depends on $\alpha$ and $\beta$ and is a continuous function of $(\alpha,\beta)$ in this domain. 
\item If $\alpha<1$ and $\beta<1$, then $p=p(\alpha,\beta)$ is the unique fixed point of $h_{(\alpha,\beta)}$, in particular $p\in (0,1)$. 
\item If $\alpha>1$ or $\beta>1$, then $p=\pp(X_n\to +\infty)$. 
\item If $\alpha=1$ and $\beta<1$, then $p=1$ (resp. if $\alpha<1$ and $\beta = 1$, then $p=0$). 
\end{itemize}
On the other hand, if $p_{i,+}=-p_{i,-}$, for all $i\ge 1$, and $\alpha=\beta=1$, then $x_n\to 1/2$ almost surely. 
\end{theo}

\noindent Note that the picture is almost complete. But we do not know whether $x_n$ converges in general when $\alpha=\beta=1$ and whether $p=\pp(X_n\to +\infty)$ only depends on $\alpha$ and $\beta$ when $\alpha>1$ and $\beta>1$.

\section{Proof of Theorem \ref{generalLLN}}
The proof is based on the observation that the sequence $(x_n,n\ge 1)$ is solution of a stochastic algorithm.
First, since $h_n$ is nonincreasing for all $n\ge 0$, $h$ is also nonincreasing. In particular it has left and right limits in all points, and there is a unique $p\in[0,1]$, such that $h(p-)\ge p\ge h(p+)$.   
Next, by definition
$$d_{n+1}=d_n +1(T^{(n+1)}_{d_n+1}<T^{(n+1)}_{g_n-1}),$$
for all $n\ge 1$. 
Therefore
$$x_{n+1}=x_n+\frac{p-x_n}{n+3}+\frac{\epsilon_{n+1}}{n+3} + \frac{r_{n+1}}{n+3} \quad \text{for all } n\ge 1,$$
where  
\begin{eqnarray}
\label{epsilon}
\epsilon_{n+1}:=1(T^{(n+1)}_{d_n+1}<T^{(n+1)}_{d_n-n-1})-h_n(x_n),
\end{eqnarray} 
is the increment of a martingale, and $r_{n+1}:=h_n(x_n)-p$. Set now $y_n:=x_n-p$. Then 
$$y_{n+1}^2=y_n^2 - \frac{2y_n^2}{n+3} + \frac{\epsilon'_{n+1}}{n+3}+\frac{s_{n+1}}{n+3}+t_{n+1}\quad \text{for all } n\ge 1,$$
where $\epsilon'_{n+1}=2\epsilon_{n+1}y_n$, $s_{n+1}=2r_{n+1}y_n$ and $t_n=\mathcal{O}(n^{-2})$. 
In particular the $\epsilon'_n$ are still the increments of a martingale. 
We claim now that    
\begin{eqnarray}
\label{hnsn}
\limsup_{n\to \infty} s_{n+1}\le 0.
\end{eqnarray}
Indeed, let $\epsilon>0$ be fixed. First if $|x_n-p|\le \epsilon$, then $|s_{n+1}|\le 2\epsilon$, since $|h_n|\le 1$. 
Now if $x_n\ge p+\epsilon$, since $h_n$ is nonincreasing, $h_n(x_n)-p \le h_n(p+\epsilon)-p$. But $h_n(p+\epsilon) \to h(p+\epsilon)\le p$, by hypothesis. This proves that $h_n(x_n)-p\le \epsilon$, for $n$ large enough. The case $x_n\le p-\epsilon$ is similar. This proves \eqref{hnsn}.   
We deduce that $s'_n:=\max(s_n,0)$ converges to $0$, as $n\to \infty$. Define next $(z_n,n\ge 1)$ by $z_1=y_1^2$ and 
$$z_{n+1}=z_n - \frac{2z_n}{n+3} + \frac{\epsilon'_{n+1}}{n+3}+\frac{s'_{n+1}}{n+3}+t_{n+1}\quad \text{for all } n\ge 1.$$
Since $|\epsilon'_n|\le 2$, it is a well known fact that $z_n\to 0$ almost surely (see \cite{Ben} or \cite[Theorem 3.III.9 p.105]{Du}).   
But on the other hand it is immediate by induction on $n$, that $0\le y_n^2\le z_n$, for all $n\ge 1$. Thus $y_n\to 0$ almost surely, 
which concludes the proof of Theorem \ref{generalLLN}. \hfill $\square$

\section{Proof of Proposition \ref{prophnh}} 
Note that $Y_0=0$ a.s., so $h(0)=1$ and $h(1)=0$. Now fix some $x\in (0,1)$.  
For $y\in \R$, set $T(n,y):= T^{(n)}_{[(n+2)y]}/n^c$. Let $\epsilon>0$ be such that $x+\epsilon<1$ and $x-\epsilon>0$. Since by hypothesis $T_1$ and $T_{-1}$ are a.s. finite, 
there exists $T>0$, such that 
$\pp[T_1\ge T]\le \epsilon$, and $\pp[T_{-1}\ge T]\le \epsilon$.  
Since $Y^{(n)}$ converges in law to $Y$, by using Skorokhod's representation theorem, we can assume that $Y^{(n)}$ converges a.s. uniformly to $Y$ 
on $[0,T]$. In other words, if  
$$\Omega_{n,\epsilon}:= \left\{\sup_{t\le T} |Y^{(n)}_t-Y_t|<\epsilon\right\},$$ 
then we can assume that $\pp(\Omega_{n,\epsilon})\to 1$, as $n\to \infty$. In particular $\pp(\Omega_{n,\epsilon})\le \epsilon$, for $n$ large enough. 
Now, since $x+\epsilon <1$, we have 
$$\pp[\Omega_{n,\epsilon},\ T(n,x)\ge T]\le \epsilon,$$ 
and similarly with $x-1$ instead of $x$. Moreover,  
$$\Omega_{n,\epsilon}\cap \left\{T(n,x)<T(n,x-1)\le T\right\} \subset \left\{T_{x-\epsilon}<T_{x-1-\epsilon}\right\}.$$
Therefore,  
\begin{eqnarray*}
h_n(x)=\pp[T(n,x)<T(n,x-1)]&\le & 3\epsilon + \pp[\Omega_{N,\epsilon}, T(n,x)<T(n,x-1)\le T]\\
&\le & 3\epsilon + \pp[T_{x-\epsilon}<T_{x-1-\epsilon}] = 3\epsilon + h(x-\epsilon), 
\end{eqnarray*} 
for $n$ large enough. We can prove similarly that $h_n(x)\ge h(x+\epsilon) - 3\epsilon$, for $n$ large enough. 
Since $\epsilon$ can be taken arbitrarily small, and since $h$ is assumed to be continuous, this concludes the proof of the proposition. \hfill $\square$

\section{Proof of Theorem \ref{LLN}}
$\bullet$ First when $\alpha<1$ and $\beta<1$, the result follows from Theorem \ref{generalLLN}, Proposition \ref{prophnh} and Theorem \ref{Dolgo}.

\vspace{0.2cm}
$\bullet$ If now $\alpha>1$ or $\beta>1$, then Zerner's result \cite{Zer1} implies 
that the cookie random walk is transient. 
We then immediately get that for any $x\in (0,1)$, $h_n(x)\to p:=\pp[X_n\to \infty]$. 
Thus we can apply Theorem \ref{generalLLN} with $h$ the function equal to $p$ on $(0,1)$, 
and such that $h(0)=1$ and $h(1)=0$.

Note that if we are in the case $\alpha>1$ and $\beta>1$, then by Lemma 8 in \cite{Zer1} we have 
$\pp[X_n>0\ \textrm{for all }n>0]>0$ and $\pp[X_n<0\ \textrm{for all }n>0]>0$. Thus $p\in (0,1)$. 
On the other hand if $\alpha>1$ and $\beta\le 1$, then $p=1$ (resp. if $\alpha\le 1$ and $\beta>1$, then $p=0$). 

\vspace{0.2cm}
$\bullet$ We consider next the case $\alpha= 1$ and $\beta<1$ (the case $\alpha<1$ and $\beta= 1$ is entirely similar), 
which belongs to the recurrent regime. 
We have to prove that $x_n\to 1$ almost surely. 
For this we use the following monotonicity result due to Zerner: 
\begin{lem}[\cite{Zer1} Lemma 15] 
Let $\pp$ be the law of a cookie random walk evolving in an environment $(p_{i,\pm},i\ge 1)$ and $\pp'$ 
the law of a cookie random walk evolving in another environment $(q_{i,\pm},i\ge 1)$. 
Assume that $p_{i,+}\ge q_{i,+}$ and $p_{i,-}\ge q_{i,-}$, for all $i\ge 1$. Then $\pp[T_y< T_x]\ge \pp'[T_y<T_x]$, for all $x<0<y$.     
\end{lem}
\noindent Two cases may appear now. Either there exists $M$ such that $\sum_{i=1}^{M-1} p_{i,+}<1\le \sum_{i=1}^M p_{i,+}$, or for all $\epsilon>0$, 
there exists $M=M(\epsilon)$ such that 
$1-\epsilon\le \sum_{i=1}^M p_{i,+}<1$. In any case, set $q_{i,-}=p_{i,-}$, for all $i\ge 1$, and $q_{i,+}=1/2$, for $i>M$. 
In addition, in the first case, set $q_{i,+}=p_{i,+}$, 
for $i\le M-1$ and $q_M^+=p_M^+-\epsilon$, with $\epsilon<p_M^+$. In the second case, set $q_{i,+}=p_{i,+}$, for $i\le M$. 
The above lemma implies that in both cases, the process $(d_n,n\ge 0)$ 
associated to $(p_{i,\pm},i\ge 1)$ stochastically dominates the process $(d'_n,n\ge 0)$ associated to $(q_{i,\pm},i\ge 1)$. 
Thus $\liminf x_n \ge p(\alpha_\epsilon,\beta)$, for all $\epsilon>0$, where $\alpha_\epsilon=\sum_{i=1}^M q_{i,+}$. 
But $\alpha_\epsilon\to 1$, when $\epsilon\to 0$. 
Moreover we already saw in the discussion above Theorem \ref{LLN} that $p(\alpha,\beta)\to 1$, when $\alpha\to 1$. 
This proves that $x_n$ converges a.s. to $1$. 

\vspace{0.2cm} 
$\bullet$ It remains to consider the case when $p_{i,+}=-p_{i,-}$, for all $i\ge 1$, and $\alpha=\beta=1$. In this case the processes $X$ and $-X$ have the same law. 
Thus $h_n(1/2)=1/2$ for all $n\ge 1$, and as an immediate consequence $(h_n(x)-1/2)(x-1/2)\le 0$, for all $x\in [0,1]$. 
In particular, we can follow the proof of Theorem \ref{generalLLN}, and we find that $s_{n+1}\le 0$, for all $n\ge 1$. Since it was all we needed in this proof, 
we conclude as well that $x_n\to 1/2$ almost surely. This finishes the proof of Theorem \ref{LLN}.  
\hfill $\square$

\end{document}